\documentclass[8pt,a4paper]{article}

\usepackage[latin1]{inputenc}
\usepackage{fancyhdr}
\usepackage{amssymb}
\usepackage{amsmath}
\usepackage{latexsym}
\usepackage{amsthm}
\usepackage[english]{babel}
\usepackage{amsfonts}
\usepackage{graphicx}
\usepackage{algorithm}
\usepackage{algorithmic}
\usepackage{subfig}
\usepackage{mathrsfs}
\usepackage{color}
\usepackage{multirow}

\newtheorem{definition}{Definition}[section]
\newtheorem{remark}{Remark}[section]
\newtheorem{example}{Example}[section]

\begin{document}

\title{Symmetric functions for fast image retrieval with persistent homology}

\author{Alessia Angeli\footnote{Dept. of Mathematics, University of Bologna, Piazza di Porta San Donato 5, Bologna, Italy}, Massimo Ferri\footnote{Dept. of Mathematics, University of Bologna, Piazza di Porta San Donato 5, Bologna, Italy} and Ivan Tomba\footnote{R\&D Dept., CA-MI S.r.l., Via Ugo La Malfa 13, Pilastro di Langhirano, Parma, Italy}}








\maketitle

\abstract{Persistence diagrams, combining geometry and topology for an effective shape description used in pattern recognition, have already proven to be an effective tool for shape representation with respect to a certain filtering function. Comparing the persistence diagram of a query with those of a database allows automatic classification or retrieval, but unfortunately, the standard method for comparing persistence diagrams, the bottleneck distance, has a high computational cost. A possible algebraic solution to this problem is to switch to comparisons between the complex polynomials whose roots are the cornerpoints of the persistence diagrams. This strategy allows to reduce the computational cost in a significant way, thereby making persistent homology based applications suitable for large scale databases. The definition of new distances in the polynomial framework poses some interesting problems, both of theoretical and practical nature. In this paper, these questions have been addressed by considering possible transformations of the half-plane where the persistence diagrams lie onto the complex plane, and by considering a certain re-normalisation the symmetric functions associated to the polynomial roots of the resulting transformed polynomial. The encouraging numerical results, obtained in a dermatology application test, suggest that the proposed method may even improve the achievements obtained by the standard methods using persistence diagrams and the bottleneck distance.}


\section{Introduction}

Persistent Homology is a recent branch of Computational Topology, which combines geometry and topology for an effective shape description used in Pattern Recognition.

It has already proven to be an effective tool for shape representation in various applications, in particular it is a fairly popular tool for comparing, retrieving and classifying objects having a natural origin, especially medical images.
Examples of such applications are liver lesion image classification (see \cite{Adcock14}), brain data analysis (\cite{Chung09} and \cite{Bend14}) and dermatology, where persistent homology has been tested with encouraging results \cite{D'Amico04,FerriStanga2010,Stanga05,Tomba17}.

The standard tool for shape comparison is the natural pseudo-distance \cite{Donatini04}, which provides a dissimilarity measure between two topological spaces under a certain filtering function, which represents the aspect of the spaces to be compared. Due to the abstract nature of its definition the natural pseudo-distance is not directly computable. Luckily, the natural pseudo-distance can be approximated by comparing persistence diagrams (also known as bar codes), which condense the essence of the pair constituted by the topological space and the filtering function in finite sets of points in the plane \cite{Fro01}: the bottleneck distance (or matching distance) between persistence diagrams yields an optimal lower bound to the natural pseudo-distance \cite{Donatini04,D'Amico10}. However, the bottleneck distance suffers from combinatorial explosion \cite{D'Amico06}, so it may become impractical to scan a large database to produce a retrieval or an automatic classification. Approximations such as image size reductions, smart organization of the database
according to the metric, progressive application of different classifiers come to help, but do not solve the problem completely. In fact, reducing images too much may result in a significant loss of information; smart organization of the database lightens the problem of a real-time retrieval as it allows to avoid some distance computations between couples of objects, but has the drawback of requiring an a-priori computation of all relative distances in the database: this results in $O(N_{db}^2)$ bottleneck distance computations (where $N_{db}$ is the number of images stored in the database), which for persistence diagrams with thousands of cornerpoints may become too demanding from a computational point of view.

A possible solution, proposed for the first time in \cite{FerriLandi99}, is to represent a persistence diagram as the set of complex roots of a polynomial; then comparison can be performed on coefficients, or on a proper subset of the coefficients \cite{DiFabioFerri2015}. From the computational point of view, this has a double advantage: first, it allows to store only a fixed number of coefficients of the polynomials instead of the persistence diagrams, sparing memory; second, it accelerates the computation of the distances, as the distance between two coefficient vectors is cheaper than the bottleneck distance. This method is intriguing, but requires the solution of an important issue: in real situations there are a lot of points near the diagonal $\Delta = \{(u, v) \in \mathbb{R}^2 : u = v\}$, due to noise, which are less meaningful in shape representation, but with a heavy impact on polynomial coefficients. Indeed, the metric of the half-plane $\Delta^+ = \{(u, v) \in \mathbb{R}^2 : u < v\}$, where persistence diagrams lie, is not to be considered Euclidean, as the diagonal acts like a unique point/object and all points close to the diagonal are close to each other in the sense introduced by the bottleneck distance. Thus, it makes sense, before passing from the persistence diagrams to the polynomial coefficients, to apply a geometric transformation which maps the entire diagonal to the origin of the complex plane and all cornerpoints close to the diagonal close to the origin. Such transformations are of interest on their own (i.e. independently on the choice to consider polynomial coefficients), as they map the persistence diagrams onto a more standard metric space, such as the complex plane with the Euclidean distance. Some transformations have been proposed in \cite{DiFabioFerri2015}: here a new transformation is proposed and compared to them in the numerical experiments.

Another practical problem, not discussed in \cite{DiFabioFerri2015}, arises when dealing with persistence diagrams with many cornerpoints: some coefficients may become very large, causing instabilities in the distance computations and more generally in the metric of the database. A solution of this problem is proposed as well.

At last it should be remembered, from polynomial theory, that a small distance between polynomial roots implies a small distance of coefficients, but the converse is false. Thus, it is possible that two very different persistence diagrams have very similar polynomial coefficients: this is the price to be paid for attempting a computational simplification. A possible solution is to use polynomial comparison as a preprocessing phase in shape retrieval, i.e. as a very fast way of getting rid of definitely far objects, so that the bottleneck distance can be computed only for a small set of candidates. However, as \cite{DiFabioFerri2015} and this paper show, the use of polynomials appears to work well in the experiments, suggesting that in practice the unfortunate case in which dissimilar persistence diagrams give rise to similar coefficient vectors is probably quite rare.

The paper is organised as follows:
Section \ref{Preliminaries} provides the theoretical background required for understanding the paper; Section \ref{SymFWD} defines a new distance between transformed persistence diagrams, in comparison to the bottleneck distance; in Section \ref{Alg} the new algorithm is stated and its computational complexity is discussed; Section \ref{ExpRes} presents the numerical simulations conducted in the dermatology application and Section \ref{Conclusions} is dedicated to further possible developments and to the conclusions.

\section{Preliminaries}\label{Preliminaries}

In persistence, each shape is viewed as a pair $(X,f)$: a topological space $X$, and a continuous function $f:X\rightarrow\mathbb{R}$, called \textit{filtering function} (or \textit{measuring function}), to define a family of subspaces $X_{\ell}=\left\{x\in X:f(x)\leq\ell\right\}$, $\ell\in\mathbb{R}$, nested by inclusion, i.e. a filtration of $X$.

Given $u$, $v\in\mathbb{R}$, $u<v$, we consider the continuous inclusion $i^{u,v}:X_{u}\rightarrow X_{v}$. This inclusion induces a homeomorphism of homology groups $i_{*}^{u,v}: H_{k}(X_{u})\rightarrow H_{k}(X_{v})$ for every $k\in\mathbb{Z}$. The image of this application are the $k$-homology classes that are born before or at level $u$ and are still alive at the level $v$ and it is called the \textit{$k$-th persistent homology group of $(X,f)$ at $(u,v)$}. If this group is finitely generated, we will denote by $\beta_{k}^{u,v}(X,f)$ its rank.

Persistent homology groups of $(X,f)$ are often described using \textit{persistence diagrams}, i.e. multisets of points whose abscissa and ordinate are, respectively, the level at which $k$-homology classes are created and the level at which they are annihilated though the filtration (\textit{proper cornerpoints}). If a homology class does not die along the filtration, the ordinate of the corresponding point is set to $+\infty$ (\textit{cornerpoint at infinity}).

Let $\Delta^{+}=\left\{(u,v)\in\mathbb{R}^{2}:u<v\right\}$, $\Delta=\left\{(u,v)\in\mathbb{R}^{2}:u=v\right\}$, $\bar{\Delta}^{+}=\Delta^{+}\cup\Delta$. \\

\begin{definition}
\emph{Proper cornerpoint}.
\\ For every point $(u,v)\in\Delta^{+}$, let us define the number $\mu_{k}(u,v)$ as the minimum over all the positive real numbers $\epsilon$, with $u+\epsilon < v-\epsilon$, of
\begin{equation*}
\beta_{k}^{u+\epsilon,v-\epsilon}(X,f)-\beta_{k}^{u-\epsilon,v-\epsilon}(X,f)-\beta_{k}^{u+\epsilon,v+\epsilon}(X,f)+\beta_{k}^{u-\epsilon,v+\epsilon}(X,f)
\end{equation*}
The finite number $\mu_{k}(u,v)$ is called \emph{multiplicity of $(u,v)\in\Delta^{+}$} for $\beta_{k}(X,f)$ and the points $(u,v)\in\Delta^{+}$ such that $\mu_{k}(u,v)>0$ are called \emph{proper cornerpoints} for $\beta_{k}(X,f)$.
\end{definition}

\begin{definition}
\emph{Cornerpoint at infinity}.
\\ For every vertical line $r$, with equation $x=w$, let us identify $r$ with the pair $(w,+\infty)$, and define the number $\mu_{k}(r)$ as the minimum, over all the positive real numbers $\epsilon$ with $w+\epsilon < \frac{1}{\epsilon}$, of
\begin{equation*}
\beta_{k}^{w+\epsilon,\frac{1}{\epsilon}}(X,f)-\beta_{k}^{w-\epsilon,\frac{1}{\epsilon}}(X,f)
\end{equation*}
The finite number $\mu_{k}(r)$ is called \emph{multiplicity of $r$} for $\beta_{k}(X,f)$ and the lines $r$ such that $\mu_{k}(r)>0$ are called \emph{cornerlines} or \emph{cornerpoints at infinity} for $\beta_{k}(X,f)$.
\end{definition}

A cornerline is identified as a cornerpoint with infinite ordinate so it escapes our processing. Therefore, in this work, we take the cornerline to a particular cornerpoint with finite coordinates.

Given a cornerline with abscissa $w$, we substitute it with the point
\begin{equation*}
(w, \, \max\left\{v \, : \, (u, v) \, \text{is a proper cornerpoint}\right\})
\end{equation*}

\medskip

In this way we represent the cornerlines with particular proper cornerpoints so in this paper, from now on, we will simply talk about cornerpoints.

\begin{definition}
\emph{Persistence diagram}.
\\ The \emph{$k$-th persistence diagram} $\mathcal{D}_k(X,f)$ is the set of cornerpoints, each counted with its multiplicity, union the set of points of $\Delta$, counted with infinite multiplicity.
\end{definition}

Persistence diagrams are usually compared computing the so called \textit{bottleneck distance} (or \textit{matching distance}) because of the robustness of these descriptors with respect to it: roughly, small changes in a given filtering function produces just a small change in the associated persistence diagram. The bottleneck distance between two persistence diagrams measures the cost of finding a correspondence between their points. In particular, the cost of taking a point $P$ to a point $P'$ is measured as the minimum between the cost of moving one point onto the other and the cost of moving both points onto the diagonal and the matching of a cornerpoint $P$ with a point of $\Delta$ can be interpreted as the destruction of the cornerpoint $P$.

\begin{definition}
\emph{Bottleneck distance} or \emph{matching distance}.
\\ Let $\mathcal{D}_{k}$ and $\mathcal{D'}_{k}$ be two persistence diagrams with a finite number of cornerpoints, the \emph{bottleneck distance} $d_{B}(\mathcal{D}_{k},\mathcal{D'}_{k})$ is defined as
\begin{equation*}
d_{B}(\mathcal{D}_{k},\mathcal{D'}_{k})=\min_{\sigma}\max_{P\in\mathcal{D}_{k}}\hat{d}(P,\sigma(P))
\end{equation*}
where $\sigma$ varies among all the bijections between $\mathcal{D}_{k}$ and $\mathcal{D'}_{k}$ and
\begin{equation*}
\hat{d}((u,v),(u',v'))=\min\left\{\max\left\{|u-u'|,|v-v'|\right\},\max\left\{\frac{v-u}{2},\frac{v'-u'}{2}\right\}\right\}
\end{equation*}
given $(u,v)\in\mathcal{D}_{k}$ and $(u',v')\in\mathcal{D'}_{k}$.
\end{definition}

\medskip

All persistence diagrams considered in this work, related to the problem about recognition of melanocytic lesions, are $0$-persistence diagrams, so in this paper, from now on, we shall drop the $0$- prefix.

\medskip

Detailed information on persistence can be found in \cite{EdHa08,EdHa09}.

\section{Symmetric functions of warped diagrams}\label{SymFWD}

Evaluating the bottleneck distance between persistence diagrams can be computationally very expensive, making its usage not practicable on large databases. In this work we propose a new procedure taken from the procedure explained in \cite{DiFabioFerri2015} and based on the preliminary idea introduced in \cite{FerriLandi99}. The problem of comparing directly two persistence diagrams through the bottleneck distance is translated into the problem of comparing complex vectors associated with each persistence diagram through an appropriate metric between vectors. The components of these complex vectors are the values of elementary symmetric functions of transformed cornerpoints of persistence diagrams.

\subsection{Warping the plane}

Given a persistence diagram, firstly we define a transformation taking the cornerpoints to the set of complex numbers.
Secondly, we consider the elementary symmetric functions of transformed cornerpoints to construct complex vectors.

Initially we analysed the transformation previously published in \cite{DiFabioFerri2015}:
\begin{itemize}
\item[] $T:\bar{\Delta}^+\rightarrow\mathbb{C}$, \quad $T(u,v)=\frac{v-u}{2}(\cos(\alpha)-\sin(\alpha)+i(\cos(\alpha)+\sin(\alpha)))$
\end{itemize}
where $\alpha=\sqrt{u^2+v^2}$.

Successively we defined a new transformation in the following way:
\begin{itemize}
\item[] $R:\bar{\Delta}^+\rightarrow\mathbb{C}$, \quad $R(u,v)=\frac{v-u}{\sqrt{2}}(\cos(\theta)+i\sin(\theta))$
\end{itemize}
where $\theta=\pi(u+v)$.

Both $T$ and $R$ are continuous maps.

We define the multiplicity of a complex number in the range of $T$ or $R$ to be the sum of the multiplicities of the points belonging to its preimage.


These transformations warp the ``diagonal'' $\Delta$ to the origin $(0,0)$ and take the cornerpoints with distance $\bar{d}$ from $\Delta$ to points with distance $\bar{d}$ from the origin. In this way the cornerpoints near $\Delta$, principally due to noise, are all moved close to the origin. In fact the goal of these transformations is to try to reduce the impact of these cornerpoints near $\Delta$ in the sums of products which will build the elementary symmetric functions of transformed cornerpoints we are going to compare.

However, these transformations have different effects and the main difference is the following: the transformation $T$ rotates the ray on which a given cornerpoint lies, proportionally to the distance from the origin $(0,0)$ of the cornerpoint itself; the transformation $R$ takes segments orthogonal to $\Delta$ to rays coming out of the origin $(0,0)$.

\begin{remark}
Regarding the transformations $T$ and $R$, it is important to observe that they are not injective on $\Delta^+$, but injectivity is restored if they are restricted to the triangle $\mathcal{T}=\{(u,v) \in \mathbb{R}^2 : 0 \leq u < v \leq 1\}$. The restriction can be achieved by normalising all filtering functions in the range $[0,1]$. We also considered several variations of these to try to better distribute the noise. However, the use of these variants did not lead to significant or global improvements in the obtained results, so we decided to keep their simplest versions.
\end{remark}

\begin{example}
When dealing, e.g., with melanocytic lesions, we report an example of how the transformations just described take the cornerpoints to complex numbers. The graphs of Fig. \ref{ImageGraph} are obtained from image $IMD435$ of the database $PH^2$, with the filtering function $Bound$ and transformations $T$ and $R$. \\

\begin{figure}[htb]
\centering
\includegraphics[width=4.72cm]{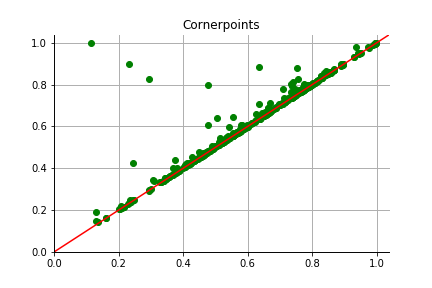}
\\
\includegraphics[width=4.72cm]{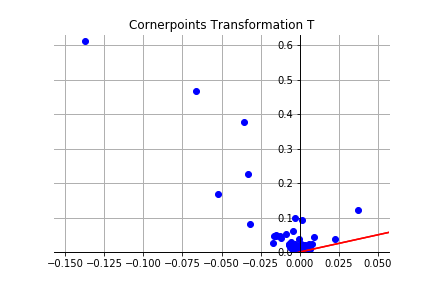}
\qquad
\includegraphics[width=4.72cm]{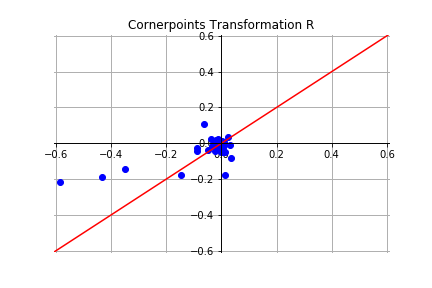}
\caption{\textit{Database $PH^{2}$, image $IMD435$, filtering function \textit{Bound}.}}
\label{ImageGraph}
\end{figure}
\end{example}

In \cite{DiFabioFerri2015} only proper cornerpoints are considered, while in this work we consider also cornerlines, because they carry information which sometimes turns essential in the present experimentation.

\begin{example}
When dealing, e.g., with melanocytic lesions, with shade of brown as a filtering function, given two images both having a fairly homogeneous hue but one light and one dark colored, these could have very similar persistence diagrams except for the abscissas of cornerlines.
\end{example}

\subsection{Comparing warped diagrams through symmetric functions}

Let $\mathcal{D}$ be a persistence diagram and $c_1=(u_1,v_1)$, ... , $c_h=(u_h,v_h)$, where $c_1$, ... , $c_h\in\mathbb{R}^2$, its cornerpoints with multiplicity $r_1$, ... , $r_h$ respectively. Let now the complex numbers $z_1$, ... , $z_h$ be obtained from $c_1$, ... , $c_h$ by one of the transformations $T$, $R$. We associate to $\mathcal{D}$ the complex vector $a_{\mathcal{D}}$ whose components are the values of elementary symmetric functions of transformed cornerpoints. So we have:
\begin{equation}
a_{\mathcal{D}}(j)=\sum_{1\leq i_{1}<i_{2}< ... <i_{j}\leq n}z_{i_{1}}z_{i_{2}} ... z_{i_{j}}, \qquad j=1,2, ... ,n
\label{FormulasVESFTC}
\end{equation}
where $n$ is the number of cornerpoints of $\mathcal{D}$ counted with multiplicity.

Now let $\mathcal{D}$ and $\mathcal{D'}$ be a persistence diagrams with $n$ and $m$ cornerpoints counted with multiplicity, respectively. If $n\neq m$ ($m<n$ say) the complex vectors $a_{\mathcal{D}}$ and $a_{\mathcal{D'}}$ corresponding respectively to $\mathcal{D}$ and $\mathcal{D'}$ have different lengths. We can manage this problem by also computing $a_{\mathcal{D'}}(j)$ with $j=m+1, ... ,n-1,n$ adding the complex number zero with multiplicity $n-m$ to the set of complex numbers obtained from the cornerpoints of $\mathcal{D'}$. In so doing, we can build two vectors of complex numbers $a_{\mathcal{D}}$ and $a_{\mathcal{D'}}$ of the same length and then we are ready to compute a distance between them.

\medskip

As for the distance to be used among complex vectors $a_{\mathcal{D}}$ and $a_{\mathcal{D'}}$, we can observe that the preliminary tests in \cite{DiFabioFerri2015} suggested that the first components of these vectors were more meaningful. Therefore we consider only the first $k$ components of these vectors, $k\in\left\{1, ... , n\right\}$. The chosen metric is the following:
\begin{equation}
d(a_{\mathcal{D}},a_{\mathcal{D'}})=\sum_{j=1}^{k}|a_{\mathcal{D}}(j)-a_{\mathcal{D'}}(j)|
\label{Distance}
\end{equation}

Applying the algorithm previously described to the persistence diagrams $\mathcal{D}$ and $\mathcal{D'}$ it happens that the difference between the orders of magnitude of $a_{\mathcal{D}}(j)$ and $a_{\mathcal{D'}}(j)$, $j\in\left\{1, ... , n\right\}$, can be relevant (more precisely, we are talking about the order of magnitude of the absolute value of the real and imaginary parts): up to six! If this occurs only for some $j\in\left\{1, ... , n\right\}$, then the distance will essentially depend only on the components with these indices, stymieing the contribution of the others.

We actually compare only the first $k$ components, but we want the distance to take all of them into account.
Moreover, if this happens for a rather big index $j$, this contrasts with our assumption of greater importance of the first components.

The main cause of this phenomenon seems to be the dependence on the number of cornerpoints, which can vary quite a lot from image to image. The summation yielding the $j$-th component will have $\binom{n}{j}$ addends but many of these addends could be null if this persistence diagram has a number of cornerpoints much smaller than $n$. Instead, if all the addends of this summation have values close to $1$, then $a_{\mathcal{D}}(j)$ will be close to $\binom{n}{j}$ with great growth at the approach of $j$ to $\frac{n}{2}$.

\begin{example}
With the database $PH^{2}$ and the filtering function \textit{LI2}, there are persistence diagrams with $47$, $77$ or $156$ cornerpoints but there are also persistence diagrams with a higher number of cornerpoints like $1444$, $1247$ or $1815$.
\end{example}

We then apply a further transformation to $a_\mathcal{D}(j)$ by extracting the $j$-th root and dividing by the number of cornerpoints of the persistence diagram considered to try to mitigate the dependence of the values of elementary symmetric function of transformed cornerpoints on that number.

More precisely, let $a_{\mathcal{D}}$ and $a_{\mathcal{D'}}$ be two complex vectors with the same length $n$ and associated respectively to persistence diagrams $\mathcal{D}$ and $\mathcal{D'}$, the defined transformation is the following:
\begin{equation}
\begin{split}
a_{\mathcal{D}}(j) \qquad &\longmapsto \qquad \tilde{a}_{\mathcal{D}}(j)= \frac{|a_{\mathcal{D}}(j)|^\frac{1}{j}}{N}e^{i arg(a_{\mathcal{D}}(j))}, \qquad j\in\left\{1, ... , k\right\} \\
a_{\mathcal{D'}}(j) \qquad &\longmapsto \qquad \tilde{a}_{\mathcal{D'}}(j)=\frac{|a_{\mathcal{D'}}(j)|^\frac{1}{j}}{M}e^{i arg(a_{\mathcal{D'}}(j))}, \qquad j\in\left\{1, ... , k\right\}
\end{split}
\label{NewTransformation}
\end{equation}
where $N$ and $M$ are respectively the number of cornerpoints of the considered persistence diagrams counted with multiplicity and $k$ is the number of components of $\tilde{a}_{\mathcal{D}}$ and $\tilde{a}_{\mathcal{D'}}$ selected for the vector comparison.

\begin{example}
In our experimentation, by applying (\ref{NewTransformation}) after one of our two transformations to the persistence diagrams, we obtain some values of complex vectors, related to the same component, with order of magnitude between $10^{-5}$ and $10^{-3}$, with prevalence of $10^{-3}$, or between $10^{-2}$ and $10^{-3}$.
\end{example}

\section{Algorithms}\label{Alg}

The algorithms below resume the principal steps of our scheme for computing the coefficients associated to the persistence diagrams extrapolated from a database of objects $db$. In \textit{Algorithm $1$} the transformation $F_{1}$ (line $2$) corresponds to one of the transformations $T$ or $R$, while in \textit{Algorithm $2$} the transformation $F_{2}$ (line $4$) corresponds to the transformation (\ref{NewTransformation}) defined above.

\medskip

\begin{tabular}{l}
\textbf{\textit{Algorithm $1$}} \\
\textbf{Input:} List $A$ of cornerpoints of a persistence diagram $\mathcal{D}$, \\
$n=\underset{\left\{A: \space A\in db\right\}}\max|A|$ \\
\textbf{Output:} List $B$ of complex numbers associated to $\mathcal{D}$\\
\begin{tabular}{r l}
\hline
{} \\
1:& \  \textbf{for each} $(u,v)\in A$ \\
2:& \   \quad \textbf{replace} $(u,v)$ by $F_{1}(u,v)$ \\
3:& \  \textbf{end for} \\
4:& \  \textbf{if} $|B|<n$ \\
5:& \   \quad \textbf{append} $n-|B|$ zeros to $B$ \\
6:& \  \textbf{end if} \\
{} \\
\hline
\end{tabular}
\end{tabular}

\medskip

\begin{tabular}{l}
\textbf{\textit{Algorithm $2$}} \\
\textbf{Input:} $n$, $B=list(z_{1}, ... , z_{n})$ associated to $\mathcal{D}$, $k\in\left\{1, ... , n\right\}$ \\
\textbf{Output:} Complex vector $\tilde{a}_{\mathcal{D}}^{k}=(\tilde{a}_{\mathcal{D}}(1), ... , \tilde{a}_{\mathcal{D}}(k))$ associated to $\mathcal{D}$ \\
\begin{tabular}{r l}
\hline
{} \\
1:& \  \textbf{set} $\tilde{a}_{\mathcal{D}}^{k}=list()$ \\
2:& \  \textbf{for} $j\in \{1,\ldots,k\}$\\
3:& \   \quad \textbf{compute} $\tilde{a}_{\mathcal{D}}=\underset{1\leq i_{1}<i_{2}< ... <i_{j}\leq n}\sum z_{i_{1}}z_{i_{2}} ... z_{i_{j}}$ \\
4:& \   \quad \textbf{replace} $\tilde{a}_{\mathcal{D}}$ by $F_{2}(\tilde{a}_{\mathcal{D}})$ \\
5:& \   \quad \textbf{append} $\tilde{a}_{\mathcal{D}}$ to $\tilde{a}_{\mathcal{D}}^{k}$ \\
6:& \  \textbf{end for} \\
{} \\
\hline
\end{tabular}
\end{tabular}

\medskip

After applying \textit{Algorithm $1$} and \textit{Algorithm $2$} to each persistence diagram belonging to the considered database, we can compute distances between the complex vectors obtained using the distance $d$ previously defined in (\ref{Distance}) and create the matrix of distances $M_{dist}$ between the objects of the database by considering $$V=\left\{\tilde{a}_{\mathcal{D}}^{k}: \tilde{a}_{\mathcal{D}}^{k} \,\mbox{complex vector associated with}\, \mathcal{D} \,\mbox{for each}\, \mathcal{D}\in db\right\}$$ and defining $M_{dist}(i,j) = d(\tilde{a}_{i}^{k},\tilde{a}_j^{k})$, for $i\in\{1, ... ,|V|\}$ and $j\in\{1, ... ,|V|\}$.

\medskip

Now let $N_{db}=|V|$ and $\mathcal{C}_{F_{1}}$ and $\mathcal{C}_{F_{2}}$ be the computational complexities of $F_{1}$ and $F_{2}$ respectively. Considering all persistence diagrams belonging to the database, the computational complexities of \textit{Algorithm $1$} and \textit{Algorithm $2$} are $\mathcal{C}_{1}=O(n\cdot\mathcal{C}_{F_{1}}\cdot N_{db})$ and $\mathcal{C}_{2}=O(k\cdot(k+n)\cdot\mathcal{C}_{F_{2}}\cdot N_{db})$ respectively.\\
The cost of \textit{Algorithm $2$} depends on how the formulas to compute the values of elementary symmetric functions of transformed cornerpoints (\ref{FormulasVESFTC}) are implemented. Following \cite{DiFabioFerri2015}, we implemented these formulas directly by computing the values of the first $k$ elementary symmetric functions of transformed cornerpoints by means of Vi\`{e}te's formulas.
Additionally, it is straightforward to verify that the cost of the transformations $F_{1}$ and $F_{2}$ is $O(1)$, so we can consider $\mathcal{C}_{1}=O(n\cdot N_{db})$ and $\mathcal{C}_{2}=O(k\cdot(k+n)\cdot N_{db})$ respectively. As regards, instead, the calculation of distances, the computational complexities related to it is $\mathcal{C}_{d}=O(k\cdot N_{db}^2)$. Therefore the whole computational cost of computing the matrix $M_{dist}$ turns out to be $\mathcal{C}=O(k\cdot(k+n)\cdot N_{db}+k\cdot N_{db}^{2})$, where $k$ is the number of components of the complex vectors used to compare them and $n$ is the maximum number of cornerpoints of persistence diagrams considered. In general, if $k$ is chosen negligible compared to the value of $n$, the cost of computing $M_{dist}$ becomes $O(n\cdot N_{db}+N_{db}^{2})$.

\medskip

Therefore, with an appropriate choice of the number $k$ of computed values of elementary symmetric functions of transformed cornerpoints, our algorithm to compare persistence diagrams results to be cheaper than using the bottleneck distance. In fact, considering two persistence diagrams $\mathcal{D}$ and $\mathcal{D'}$ with $N$ and $M$, respectively, the number of their cornerpoints, the cost of computing the bottleneck distance $d_{B}(\mathcal{D},\mathcal{D'})$ is $\mathcal{C}_{B}=O(\max(N,M)^\frac{3}{2}\cdot log(\max(N,M)))$ \cite{Kerber16}. Thus, in the case of the Bottleneck distance, the cost of computing $M_{dist}$ is $O(n^\frac{3}{2}\cdot \log(n)\cdot N_{db}^2)$.

\medskip

Moreover, consider the case a new object is to be compared with those of the database. This situation may occur both in the case the user wants to retrieve the first nearest neighbors of the new object, or in the case the new object needs to be classified on the base of its relative distances with respect to the objects of the database. This requires an almost instant computation of $N_{db}$ distances to provide the user with a real-time answer, thus in these applications the algorithm based on the polynomial coefficients may be more suitable for this reason as well.

\section{Experimental results}\label{ExpRes}

In this work we used the dataset $PH^{2}$ \cite{Medonca13}, a public database available online since the EMBC $2013$ conference. The $PH^{2}$ database contains a total number of $200$ melanocytic lesions, including $80$ common nevi, $80$ atypical nevi, and $40$ melanomas. The dermoscopic images were obtained under the same conditions using a magnification of $20\times$. They are $8$-bit RGB color images with a resolution of $768\times 560$ pixels. In particular we tested the algorithm previously described considering the images of the database $PH^{2}$ and $19$ features, $11$ of which filtering functions obtained using $1$-dimensional \textit{persistent homology} theory, and $8$ of them related to the \textit{ABCDE} analysis.

The features considered are summarised in Table \ref{Features}.

\begin{table}[htb]\scriptsize
\centering
\begin{tabular}{cc}
\hline
\textbf{persistence features}	& \textbf{ABCDE features} \\
\hline
Light Intensity (\textit{LI})            & Colour Histogram (\textit{Histo}) \\
Blue (\textit{B})		                     & Form Factor (\textit{FF}) \\
Green (\textit{G})		                   & Haralick's Circularity (\textit{HC}) \\
Red (\textit{R})		                     & Asymmetry (\textit{Asym}) \\
Excess Blue (\textit{ExcB})		           & Ellipticity (\textit{Elt}) \\
Excess Green (\textit{ExcG})	           & Eccentricity (\textit{Ecc}) \\
Excess Red (\textit{ExcR})	             & Diameter (\textit{Diam}) \\
Light Intensity 2 (\textit{LI2})		     & Colour Entropy (\textit{Entr}) \\
Boundary Light Intensity (\textit{BLI})  & \\
Boundary (\textit{Bound})		             & \\
Boundary 2 (\textit{Bound2})	           & \\
\hline
{}
\end{tabular}
\caption{Features.}
\label{Features}
\end{table}

The procedure used to compute the persistence diagrams from a dermoscopic image and the ABCDE features are the same as described in \cite{Tomba17}. Of the persistent homology related features, the filtering functions \textit{B}, \textit{G}, \textit{R}, \textit{ExcB}, \textit{ExcG}, \textit{ExcR} have already been described in \cite{Tomba17}. Here, we recall that the persistence related features, i.e. the 1-dimensional persistence diagrams, are obtained from a segmented image as follows:
\begin{itemize}
  \item computation of the global graph, whose nodes are all pixels of the segmentation mask, with the 6-connectivity notion;
  \item definition of the filtering functions on the nodes of the graphs;
  \item computation of the 1-dimensional PBN functions through the algorithm described in \cite{Tomba17} and the references therein.
\end{itemize}
We also recall that, for computational reasons, prior to these operations the image and the corresponding segmentation mask needed to be rescaled into smaller size images by taking averages onto squared blocks of pixels. Since in \cite{Tomba17} the use blocks of $4 \times 4$ pixels, provided good classification results for images of a similar size to those of the $PH^{2}$ dataset, here the same choice has been made. The global graph has also been used to compute the light intensity features, which are obtained from the R, G, B filtering functions as follows:
\begin{itemize}
\item \textit{LI} = (\textit{R} + \textit{G} + \textit{B})/3;
\item \textit{LI2} = 255 - \textit{LI}.
\end{itemize}

The \textit{LI2} filtering function represents the inverse filtration of the \textit{LI} filtering function. In a similar way, \textit{Bound2} is obtained by taking the inverse filtration of the distance from the barycenter on the boundary graph, described in \cite{Tomba17}. Since the boundary graph has much fewer nodes than the boundary graph, no subdivision into blocks was needed. The \textit{BLI} feature is obtained by considering the light intensity on the subgraph of the global graph composed by the pixels whose distance with respect to 4-connectivity is lower or equal to a certain constant C. This is obtained by taking one fifth of the of the maximal distance of the global graph pixels from the boundary graph. This choice allows to consider smaller blocks with respect to the global graph and to get a more detailed analysis of the texture of the skin lesion near the boundary: in this case, $2 \times 2$ block sizes have been used.
Thus, for each method considered, 19 distance matrices have been computed, each containing the relative distances between the images of the database (of course, the main diagonal of these matrices is null):
\begin{itemize}
  \item 11 matrices for the 1-dimensional filtering functions described above, using 1-dimensional
matching distances between the corresponding PBNs;
  \item 8 matrices for the ABCDE features, using the same distances as in \cite{Tomba17}.
\end{itemize}

Three different methods have been considered for the computation of the 11 persistent homology distance matrices:
\begin{itemize}
  \item[\textbf{M1:}] Bottleneck distance with persistence diagrams;
  \item[\textbf{M2:}] Algorithm 1 with transformation $T$ and \textit{Algorithm $2$} with $k=10$;
  \item[\textbf{M3:}] Algorithm 1 with transformation $R$ and \textit{Algorithm $2$} with $k=10$.
\end{itemize}

Re-ordering the $j$-th column of the distance matrices, and considering only the first nearest neighbors, it is possible to compute the retrieval of the $j$-th image of the database. From the retrieval, using the histopathological diagnosis of the images provided by the $PH^2$ database, it is possible to compute an automatical diagnosis and to compare it with the effective diagnosis of the query, therewith obtaining specificity and sensibility results of each method, with respect to a specific filtering function, on the entire database.

By considering linear combinations with non-negative coefficients of the distance matrices, it is possible to obtain global distance matrices, which in turn provide global retrievals and global sensitivity and specificity results.

All results will be reported in percentage terms of:
\begin{itemize}
\item \textit{accuracy}: number of images classified correctly on the total images of the database;
\item \textit{sensitivity}: number of melanomas classified correctly on the total number of melanomas;
\item \textit{specificity}: number of nevi classified correctly on the total number of nevi.
\end{itemize}

\subsection{Test 1}

Initially we report in Table \ref{ResultsTest1} the results related of each individual filtering function about test carried out considering database $PH^{2}$ and the three methods \textbf{M1}, \textbf{M2} and \textbf{M3}.

\begin{table}[htb]\scriptsize
\begin{center}
\begin{tabular}{|c|c|c|c|c|c|c|c|c|c|}
\cline{5-10}
\multicolumn{4}{c|}{} & \multicolumn{6}{c|}{k=10} \\
\cline{2-10}
\multicolumn{1}{c|}{} & \multicolumn{3}{c|}{\textit{bottleneck}} & \multicolumn{3}{c|}{\textit{transformation $T$}} & \multicolumn{3}{c|}{\textit{transformation $R$}} \\
\cline{2-10}
\multicolumn{1}{c|}{} & acc. & sens. & spec. & acc. & sens. & spec. & acc. & sens. & spec. \\
\hline
LI     & 73.5 & 62.5 & 76.3 & 77.5 & 80.0 & 76.9 & 72.5 & 60.0 & 75.6 \\
B      & 75.5 & 52.5 & 81.3 & 76.5 & 70.0 & 78.1 & 72.0 & 52.5 & 76.9 \\
G      & 74.5 & 60.0 & 78.1 & 73.5 & 72.5 & 73.8 & 71.0 & 70.0 & 71.3 \\
R      & 79.5 & 75.0 & 80.6 & 74.5 & 77.5 & 73.8 & 72.5 & 65.0 & 74.4 \\
ExcB   & 91.5 & 87.5 & 92.5 & 91.0 & 95.0 & 90.0 & 91.5 & 92.5 & 91.3 \\
ExcG   & 79.5 & 90.0 & 76.9 & 81.5 & 87.5 & 80.0 & 80.0 & 87.5 & 78.1 \\
ExcR   & 84.0 & 87.5 & 83.1 & 89.5 & 92.5 & 88.8 & 88.0 & 92.5 & 86.9 \\
LI2    & 82.5 & 72.5 & 85.0 & 80.5 & 70.0 & 83.1 & 84.5 & 60.0 & 90.6 \\
BLI    & 79.5 & 65.0 & 83.1 & 82.0 & 75.0 & 83.8 & 79.0 & 67.5 & 81.9 \\
Bound  & 84.5 & 80.0 & 85.6 & 77.5 & 87.5 & 75.0 & 76.5 & 85.0 & 74.4 \\
Buond2 & 82.0 & 67.5 & 85.6 & 77.5 & 85.0 & 75.6 & 74.0 & 87.5 & 70.6 \\
\hline
\multicolumn{1}{c}{} \\
\end{tabular}
\caption{\textit{Filtering function results, database $PH^{2}$.}}
\label{ResultsTest1}
\end{center}
\end{table}



\begin{remark}
In Table \ref{ResultsTest1} we reported only the singular results about \textbf{persistence features} because the new algorithm does not concern the procedure by which the results relative to the individual \textbf{ABCDE features} are computed.
\end{remark}


Observing the results shown in Table \ref{ResultsTest1} we can affirm that in this case the results obtained applying the new algorithm with $k=10$, considering both transformation $T$ and transformation $R$, are completely comparable to those obtained by using the bottleneck distance. Furthermore, among the methods used, there does not seem to be one that globally gives significantly better results than the others.

\subsection{Test 2}

In the previous tests we examined only the results related of each individual filtering function while now we are going to consider the set of all the features reported in Table \ref{Features}: \textbf{persistence features} (filtering functions) and \textbf{ABCDE features}. In particular we report in Table \ref{ResultsTest2} the global results obtained considering database $PH^{2}$ and applying the new algorithm varying or the value of $k$ such that $k\in\left\{5,10,20,50\right\}$ or the transformation to take cornerpoints to complex number among $T$ and $R$. Then we compare these results to each other and to those obtained by applying the bottleneck distance.

\medskip

The global results reported below are obtained by weighted averages of the distances relative to the different filtering functions and features. The values attributed to these weights may be due to theoretical or experimental reasons and may also be changed during construction using a relevance feedback. In particular in this case the weights to be assigned to the different features are computed with a percentage optimization on all the $19$ features considered.

\begin{table}[htb]\scriptsize
\begin{center}
\begin{tabular}{|c|c|c|c|c|}
\cline{3-5}
\multicolumn{2}{c|}{} & acc. & sens. & spec.\\
\hline
\multicolumn{2}{|c|}{\textit{bottleneck}} & 96.0 & 95.0 & 96.3 \\
\hline
\multirow{4}*{\textit{transformation $T$}}
& k=5 & 96.0 & 90.0 & 97.5 \\
\cline{2-5}
& k=10 & 95.5 & 92.5 & 96.3 \\
\cline{2-5}
& k=20 & 97.5 & 95.0 & 98.1 \\
\cline{2-5}
& k=50 & 97.0 & 92.5 & 98.1 \\
\hline
\multirow{4}*{\textit{transformation $R$}}
& k=5 & 97.5 & 95.0 & 98.1 \\
\cline{2-5}
& k=10 & 97.5 & 97.5 & 97.5 \\
\cline{2-5}
& k=20 & 97.0 & 92.5 & 98.1 \\
\cline{2-5}
& k=50 & 97.0 & 95.0 & 97.5 \\
\hline
\multicolumn{1}{c}{} \\
\end{tabular}
\medskip
\caption{\textit{Global results, database $PH^{2}$.}}
\label{ResultsTest2}
\end{center}
\end{table}

Observing the results shown in Table \ref{ResultsTest2} we can affirm that in this case also the global results obtained applying the new algorithm with $k\in\left\{5,10,20,50\right\}$, both considering transformation $T$ and transformation $R$, are completely comparable to those obtained by using the bottleneck distance. Furthermore, there does not seem to be a significant difference in the results obtained when $k$ changes in $\left\{5,10,20,50\right\}$ therefore it is advisable to consider a value of $k$ small to avoid raising the execution time of the algorithm too much.

Unfortunately, we could not test the algorithms of the considered methods on the same platform, so the computational times are not directly comparable. Nevertheless, it is significant that the computational time required for computing a distance matrix with method $\textbf{M1}$ was approximately 10 times higher than that required by the methods $\textbf{M2}$ and $\textbf{M3}$.

\section{Further developments and conclusions}\label{Conclusions}

The strategy of using symmetric functions of warped persistence diagrams is an interesting way to solve the problem of the high computational cost of the bottleneck distance, keeping the persistent homology framework intact and allowing scalability on large databases.

The warping plane transformations are of interest on their own, as they map the noisy cornerpoints near the diagonal to a neighbourhood of the origin, simplifying the geometry of the plane where cornerpoints lie: for example, we believe that defining an appropriate matching distance between the transformed diagrams could lead to interesting results as well, and we propose to investigate this in future researches.

The symmetric functions, unfortunately, appear to be highly dependent on the number of cornerpoints, thus the implementation of this method requires to filter the noise in some way. The solution proposed in Section \ref{SymFWD}, based on taking the $j$-th square root of the $j$-th symmetric function, works well in the experiments, but requires further investigation as it only seems to attenuate the problem.

Overall, the proposed method may serve to fill the gap between simulation and real practice on large scale databases for persistent homology based applications, but requires to be validated by further experiments and possibly improved by theoretical adjustments.

\bibliography{mybibfile}%








\end{document}